\documentclass{conm-p-l}

\def\NN{{\hbox{{\rm I$\!$\rm N}}}}

\newtheorem{theorem}{Theorem}[section]
\newtheorem{lemma}[theorem]{Lemma}

\numberwithin{equation}{section}

\begin{document}

\title{Symmetric Products and $Q$--manifolds}
Ä
\author{Alejandro Illanes}
\address{Instituto de Matem\'aticas,
Circuito Exterior,
Ciudad Universitaria,
M\'exico, D.F., C.P. 04510,
M\'exico}
\email{illanes@gauss.matem.unam.mx}
\thanks{The research of this paper was done
with the grant of U.~N.~A.~M., P.~A.~P.~I.~I.~T., IN 103095}

\author{Sergio Mac\'{\i}as}
\address{Instituto de Matem\'aticas,
Circuito Exterior,
Ciudad Universitaria,
M\'exico, D.F., C.P. 04510,
M\'exico}
\email{macias@servidor.unam.mx}

\author{Sam B.~Nadler, Jr.}
\address{Department of Mathematics,
West Virginia University,
P.O. Box 6310, Morgantown,
WV 26506--6310, U.~S.~A.}
\email{nadler@math.wvu.edu}



\begin{abstract}
An example is given of a compact absolute retract that is
not a Hilbert cube manifold but whose second symmetric product is the
Hilbert cube. A factor theorem is given for $n^{th}$ symmetric product of
the
cartesian product of any absolute neighborhood retract with the Hilbert
cube. A short proof is included of the known result that symmetric
products
preserve the property of being a compact Hilbert cube manifold (the
theorem
is proved here for all Hilbert cube manifolds).
\end{abstract}

\maketitle

\section{Introduction}
 Let $X$ be metric space, and let $n$ be a positive
integer. Let
$F_n(X)=\{A\subset X\ |\ A \hbox{\rm\ has at most $n$ points}\}$ with the
Vietoris
topology (which is the same as the Hausdorff metric topology)--see [11].
Then $F_n(X)$ is 
the {\it $n^{th}$ symmetric product of $X$.} Symmetric
products were first
defined and investigated by Borsuk and Ulam [3].

We let $Q$ denote the Hilbert cube. It is known that if $M$ is a compact
$Q$--manifold,
then
$F_n(M)$ is a $Q$--manifold ([5]; we present a short proof for all
$Q$--manifolds in
section 2). We show that the converse is false for
$n=2$; in particular, we give an example of a compact absolute retract $X$
that is not a $Q$--manifold but for which $F_2(X)$ is the Hilbert cube. We
also
prove
a factor theorem
for $F_n(X \times Q)$ when $X$ is any absolute neighborhood retract.

Our notation is standard. Nevertheless, we note the following: $\NN$
denotes the
set of positive
integers; $\mathcal{I}$ denotes the closed unit interval $[0,1]$; $AR$
($ANR$)
stands for
absolute (neighborhood) retract without assuming compactness; $1_S$
denotes the
identity map on a space $S$;
$\approx$ means ``is homeomorphic to".

Let $(Y,d)$ be a compact metric space. A closed subset, $A$, of $Y$ is a
{\it $Z$--set in $Y$}
provided that for each $\varepsilon>0$, there is a map
$f_\varepsilon\colon
Y\to Y\setminus
A$ such that $d(y,f_\varepsilon(y))<\varepsilon$ for each $y\in Y$ [4]. A
{\it $Z$--map} is a
map whose image is a $Z$--set in the range space.

Given a map $f\colon X\to Y$ and a positive integer $n$, we denote the
{\it
induced map} from
$F_n(X)$ into $F_n(Y)$ by $F_n(f)$; in other words, $(F_n(f))(A)=f(A)$ for
each $A\in F_n(X)$.
\section{Short Proof of  $Q$--manifold Theorem} 
We give a short,
self--contained proof of
the fact that symmetric products preserve $Q$--manifolds. The theorem is
in [5] for the
compact case with a somewhat different proof.

The following lemma is easily proved using the induced maps $F_n(f)$
defined above.

\begin{lemma}If $Y$ is a $Z$--set of the compactum $X$ then for each
$n\in\NN$,
$F_n(Y)$ is a $Z$--set in $F_n(X)$.
\end{lemma} 

The next lemma is easy to prove using Lemma 1.

\begin{lemma} Let $X$ be a compactum. Fix $n\in\NN$. If $1_X$ is the
uniform limit of a sequence $\{f_k\}_{k=1}^\infty$ of $Z$--maps, then
$1_{F_n(X)}$ is the
uniform limit of the sequence $\{F_n(f_k)\}_{k=1}^\infty$ of $Z$--maps.
\end{lemma}

\begin{lemma}If $X$ is a contractible continuum, then $F_n(X)$ is
contractible for each
$n\in\NN$.
\end{lemma}

\begin{proof} Fix $x_0\in X$ and $n\in\NN$. Let $G\colon
X\times\mathcal{I}\to X$ be a
homotopy contracting $X$ to $x_0$. Then define $\mathcal{G}\colon
F_n(X)\times\mathcal{I}\to
F_n(X)$ as follows:
$\mathcal{G}(A,t)=G(A\times\{t\})$ for each $(A,t)\in
F_n(X)\times\mathcal{I}$. Clearly, $\mathcal{G}$
contracts $F_n(X)$ to the point $\{x_0\}$.  
\end{proof} 

\begin{theorem} {\rm [5, p.~223]} For each $n\in\NN$, $F_n(Q) \approx Q$.
\end{theorem}

\begin{proof}  Assume that $Q=\prod_{j=1}^\infty\mathcal{I}_j$. Fix
$n\in\NN$. Let
$\varepsilon>0$ and let $k\in\NN$ be such that
${1\over 2^k}<\varepsilon$. Define the map $G\colon
F_n(Q)\times\mathcal{I}\to
F_n(Q)$
by

\begin{equation}
\begin{split}
G((\{(x_\ell^1)_{\ell=1}^\infty,\ldots,(x_\ell^n)_{\ell=1}^\infty\},t))\!=&
\{(x_1^1,\ldots,x_{k-1}^1,x_k^1,(1-t)x_{k+1}^1,(1-t)x_{k+2}^1,\ldots),\\
\ldots,&(x_1^n,\ldots,
x_{k-1}^n,x_k^n,(1-t)x_{k+1},(1-t) x_{k+2}^n,\ldots)\}.
\end{split}
\end{equation}

Observe that for each $A\in F_n(Q)$, we have that $G(A,0)=A$ and
diam$(G(\{A\}\times\mathcal{I}))<{1\over 2^k}<\varepsilon$; also,
$G(F_n(Q)\times\{1\})$ is
homeomorphic to $F_n(\mathcal{I}^k)$. Thus, since $F_n(\mathcal{I}^k)$ is
an $AR$ (see [6],
Korollar 2), we have
that $F_n(Q)$ is an $ANR$ (see [9], Lemma 1). Therefore, since
$F_n(Q)$ is contractible by Lemma 3, $F_n(Q)$ is an $AR$ (see [8],
Theorem 7.1).

Observe that for each $m\in\NN$, the map $f_m\colon Q \to Q$ given by
$f_m((x_\ell)_{\ell=1}^\infty)=(x_1,\ldots,x_m,0,0,\ldots)$ is a $Z$--map
and the
sequence $\{f_m\}_{m=1}^\infty$ converges uniformly to the identity map
$1_Q$ of
$Q$. Hence, by Lemma $2$, $\{F_n(f_m)\}_{m=1}^\infty$ is a sequence of
$Z$--maps that converges uniformly to $1_{F_n(Q)}$. Therefore, by
Theorem 1 of [12] and
22.1 of [4], we have that $F_n(Q)\approx Q$.  
\end{proof} 

The following theorem is proved in [5, p.~222] for compact $Q$-manifolds:  
 
\begin{theorem} If $M$ is a $Q$--manifold, then
$F_n(M)$ is a
$Q$--manifold for each $n\in\NN$.
\end{theorem}

\begin{proof}  Let $A\in F_n(M)$. By 37.2 of [4], $M\approx P\times Q$ by
a
homeomorphism
$h$, where $P$ is a polyhedron (not necessarily compact). Note that there
is a compact $AR$,
$B$, in
$P$ that is a neighborhood of the projection of $h(A)$ into $P$. Now, by
22.1 and 44.1 of [4],
$B\times Q$ is a Hilbert cube neighborhood in $P\times Q$ of $h(A)$.
Therefore,
$F_n(h^{-1}(B\times Q))$ is a neighborhood of $A$ in $F_n(M)$ and, by
Theorem 2.4,
$F_n(h^{-1}(B\times Q))\approx Q$.  
\end{proof} 

\section{Example} We present an example of a compact $AR$, $X$, such that
$X$ is
not a $Q$--manifold but $F_2(X)$ is the Hilbert cube.

Let $X=Q_1\cup Q_2$, where $Q_1\approx Q$, $Q_2\approx Q$, and
$Q_1\cap Q_2=\{p\}$. Since $p$ is a cut point of $X$, $X$ is not a
$Q$--manifold. We
show that $F_2(X)\approx Q$. To this end, let

\begin{equation}
\mathcal{K} = \{\{x_1, x_2\} \in F_2 (X): x_1 \in Q_1  {\rm ~ and ~ } x_2
\in Q_2\}
\end{equation}

\noindent and note that 

\begin{equation}
F_2(X)=F_2(Q_1)\cup F_2(Q_2)\cup  \mathcal{K}.
\end{equation}

Let $f\colon Q_1\times Q_2\to\mathcal{K}$ be given by
$f((x_1,x_2))=\{x_1,x_2\}$; we
see that $f$ is
a homeomorphism of $Q_1\times Q_2$ onto $\mathcal{K}$. Therefore,
$\mathcal{K}\approx Q$.

Next, observe that 

\begin{equation}
F_2(Q_1)\cap\mathcal{K}=\{\{x,p\}\in F_2(X) : x\in Q_1\}.
\end{equation}

\noindent Define $g:F_2(Q_1) \cap \mathcal{K} \rightarrow$ $Q_1$ as
follows: 
$g(\{x,p\}) = x$ if $x \neq p$, and $g(\{p\}) = p$.  Then, $g$ is a
homeomorphism of $F_2(Q_1) \cap \mathcal{K}$ onto $Q_1$.  Therefore,
$F_2(Q_1) \cap
\mathcal{K} \approx$ $Q$.

Now, we show that
$F_2(Q_1)\cap\mathcal{K}$ is a $Z$--set in $F_2(Q_1)$. Fix
$\varepsilon>0$.
Let $t\in (0,1)$ be such that $1-t<\varepsilon$. Without loss of
generality, we
assume that $Q_1$ is the standard Hilbert cube (i. e.,
$\prod_{k=1}^\infty\mathcal{I}_k$)
and that $p=(1,1,\ldots)$.  Let
$h_t\colon F_2(Q_1)\to F_2(Q_1)$ be given by $h_t(\{x,x'\})=\{t\cdot
x,t\cdot x'\}$. It is
easy to see that $h_t$ is continuous. Furthermore,
$h_t(F_2(Q_1))\cap(F_2(Q_1)\cap\mathcal{K})=\emptyset$ and $h_t$ moves
points
less than $\varepsilon$. Thus, $F_2(Q_1)\cap\mathcal{K}$ is a
$Z$--set in $F_2(Q_1)$.

By Theorem 2.4, $F_2(Q_1)\approx Q$. Also, we have shown that
$\mathcal{K}\approx Q$ and
that $F_2(Q_1)\cap\mathcal{K}$ is a $Z$--set in $F_2(Q_1)$. Therefore,
$F_2(Q_1)\cup\mathcal{K}\approx Q$ (Theorem 1 of [7]).

Now, observe that

\begin{equation}
\left[F_2(Q_1)\cup\mathcal{K}\right]\cap F_2(Q_2)=\left[F_2(Q_1)\cap
F_2(Q_2)
\right]\cup\left[\mathcal{K}\cap F_2(Q_2)\right]=\mathcal{K}\cap
F_2(Q_2).
\end{equation}

\noindent As before, it can be shown that $\mathcal{K}\cap F_2(Q_2)\approx
Q$
and that $\mathcal{K}\cap F_2(Q_2)$  is a $Z$--set in $F_2(Q_2)$. Thus,
since
we have shown
that $F_2(Q_1)\cup\mathcal{K}\approx Q$ and since $F_2(Q_2)\approx Q$ (by
Theorem 2.4),
we see that $\left[F_2(Q_1)\cup \mathcal{K}
\right]\cup F_2(Q_2)\approx Q$ (Theorem 1 of [7]). Therefore, since
$\left[F_2(Q_1)\cup\mathcal{K}\right]\cup F_2(Q_2)=F_2(X)$, we have proved
that
$F_2(X)\approx Q$.\par

\section{Factor Theorem} We use the following lemma in the proof of the
factor theorem
(Theorem 4.2). The symbol $co$ denotes the convex hull operator. We note
following fact: Let
$\mathcal{B}$ be strictly convex Banach space and let $p\in\mathcal{B}$;
then, each nonempty
compact convex
subset of $\mathcal{B}$ has a unique point nearest $p$ (Lemma 2.1 of [1]).

\begin{lemma} Let $Q$ be a convex Hilbert cube in the
Hilbert space $\ell_2$,
the metric for $Q$ being obtained from the strictly convex norm on
$\ell_2$.
Fix $p\in Q$. For each $K\in F_n(Q)$, let $\eta_p(K)$ be unique point
of $co(K)$
nearest $p$. Then, $\eta_p\colon F_n(Q)\to Q$ is continuous;
furthermore, if
$\{q_1,\ldots,q_r\}\in F_n(Q)$, $\eta_p(\{q_1,\ldots,q_r\})=q$, and
$t\in\mathcal{I}$, then

\begin{equation}
\eta_p(\{(1-t)q_1+tq,\ldots,(1-t)q_r+tq\})=q.
\end{equation}
\end{lemma}

\begin{proof}  The continuity of $\eta_p$ is well known. To prove the
second
part of the lemma,
let $E=\{q_1,\ldots,q_r\}$ and let $F=\{(1-t)q_1+tq,\ldots,(1-t)q_r+tq\}$.
Since $q\in
E$, clearly,
$F\subset co(E)$. Hence, $co(F)\subset co(E)$. Also, $q\in co(F)$, which
is
seen as follows:

Since $q\in E$, there are
$s_1,\ldots,s_r\in\mathcal{I}$ such that $\sum_{j=1}^rs_j=1$ and
$q=\sum_{j=1}^rs_jq_j$; then,
\begin{equation}
\sum_{j=1}^rs_j[(1-t)q_j+tq]=(1-t)\sum_{j=1}^rs_jq_j
+tq\sum_{j=1}^rs_j=(1-t)q+tq=q;
\end{equation}
hence, $q\in co(F)$. Therefore, since $co(F)\subset co(E)$ and
$\eta_p(E)=q$, it is clear that
$\eta_p(F)=q$.  
\end{proof} 

\begin{theorem} If $X$ is an $ANR$, then $F_n(X\times Q)\approx
F_n(X)\times Q$
for each $n\in\NN$.
\end{theorem}

\begin{proof}  Fix $n\in\NN$. Since $X\times Q$ is a $Q$--manifold (44.1
of [4]),
$F_n(X\times Q)$ is a
$Q$--manifold by Theorem~2. Note that $F_n(X)$ is an $ANR$ (by [10]
using 10.1 and 10.2
of [2], pp. 96--97); hence, $F_n(X)\times Q$ is a $Q$--manifold (44.1
of [4]).
Therefore, in order to show that $F_n(X\times Q)\approx
F_n(X)\times Q$, it suffices to
obtain a $CE$ map, $f$, from $F_n(X\times Q)$ onto $F_n(X)\times Q$
(43.2 of [4]).

Assume that $Q$, $p$ and $\eta_p$ are as in Lemma~{4.1}. Define $f\colon
F_n(X\times Q)\to
F_n(X)\times Q$ as follows: For each
$\{(x_1,q_1),\ldots,(x_\ell,q_\ell)\}\in
F_n(X\times Q)$,

\begin{equation}
f\left(\{(x_1,q_1),\ldots,(x_\ell,q_\ell)\}\right)=\left(\{x_1,\ldots,x_\ell\},
\eta_p(\{q_1,\ldots,q_\ell\})\right)
\end{equation}

It follows easily that $f$ is continuous. Also, $f$ maps onto
$F_n(X)\times Q$ since for
each $(\{x_1,\ldots,x_j\},q)\in F_n(X)\times Q$,
$f(\{(x_1,q),\ldots(x_j,q)\})=(\{x_1,\ldots,x_j\},q)$.

We prove that $f$ is a $CE$ map by showing that the fibers of $f$ are
contractible ([4], p.~91).
Fix $A\in F_n(X)\times Q$, where 
$A=(\{x_1,\ldots,x_\ell\},q)$. For any
$\{(y_1,q_1),\ldots,(y_r,q_r)\}\in f^{-1}(A)$ and any $t\in\mathcal{I}$,
let
\begin{equation}
G\left((\{(y_1,q_1),
\ldots,(y_r,q_r)\},t)\right)=\{(y_1, (1-t)q_1+tq),\ldots,(y_r,
(1-t)q_r+tq)\}.
\end{equation}

We show that $G$ maps $f^{-1}(A) \times \mathcal{I}$ into $f^{-1}(A)$.
Let $\{(y_1,
q_1), \ldots, (y_r, q_r)\} \in f^{-1}(A)$ and let $t \in \mathcal{I}$.
Since $f(\{(y_1,q_1),\ldots,$ $(y_r,q_r)\})=A$,
$\{y_1,\ldots, y_r\}= \{x_1,\ldots,x_\ell\}$ and
$\eta_p(\{q_1,\ldots,q_r\})=q$; hence, by Lemma~{4.1},
$\eta_p(\{(1-t)q_1+tq,\ldots,(1-t)q_r+tq\})=q$.  
Therefore, $G\left((\{(y_1,q_1),\ldots,(y_r,q_r)\},t)\right)$ $\in
f^{-1}(A)$. 

Also, note that for
any $\{(y_1,q_1),\ldots,(y_r,q_r)\}\in f^{-1}(A)$,
\begin{equation}
G\left((\{(y_1,q_1),\ldots,(y_r,q_r)\},1)\right)=\{(y_1,q),\ldots,(y_r,q)\}=
\{(x_1,q),\ldots,(x_\ell,q)\},
\end{equation}
\noindent the second equality being due to the fact that $f(\{(y_1, q_1),
\ldots, (y_r, q_r)\}) = A$.

Hence, we have shown that
$G$ is a homotopy contracting $f^{-1}(A)$ to the point $\{(x_1,q),\ldots,$
$(x_\ell,q)\}$ in $f^{-1}(A)$.
Therefore, we have proved that $f$ is a $CE$ map.  
\end{proof}

\end{document}